\documentclass{article}

\usepackage{amsmath,amsfonts,bm}

\def\eqref#1{equation~\ref{#1}}

\def\1{\bm{1}}

\DeclareMathAlphabet{\mathsfit}{\encodingdefault}{\sfdefault}{m}{sl}
\SetMathAlphabet{\mathsfit}{bold}{\encodingdefault}{\sfdefault}{bx}{n}

\newcommand{\hf}{{\frac 12}}

\newcommand{\grad}{\ensuremath{\nabla}}

\newcommand{\bfA}{{\bf A}}

\newcommand{\bfG}{{\bf G}}

\newcommand{\bfI}{{\bf I}}

\newcommand{\bfK}{{\bf K}}

\newcommand{\bfU}{{\bf U}}

\newcommand{\bfW}{{\bf W}}
\newcommand{\bfX}{{\bf X}}

\newcommand{\bfa}{{\bf a}}
\newcommand{\bfb}{{\bf b}}

\newcommand{\bfx}{{\bf x}}
\newcommand{\bfy}{{\bf y}}

\newcommand{\bfr}{{\bf r}}
\newcommand{\bff}{{\bf f}}

\newcommand{\bfw}{{\bf w}}
\newcommand{\bfz}{{\bf z}}

\newcommand{\bftheta}{{\boldsymbol \theta}}
\newcommand{\bfTheta}{{\boldsymbol \Theta}}

\newcommand{\norm}[1]{\ensuremath{\left\|#1\right\|}}

\usepackage{amsmath}
\usepackage{amsthm}
\newtheorem*{kelley-theorem*}{Theorem 5.7 - Kelley, 2011}

\usepackage{microtype}
\usepackage{graphicx}
\usepackage{subfigure}
\usepackage{booktabs} 

\usepackage{hyperref}

\usepackage[accepted]{icml2021}

\icmltitlerunning{Deep Neural Network Accelerated Implicit Filtering}

\begin{document}

\twocolumn[
\icmltitle{Deep Neural Network Accelerated Implicit Filtering}



\icmlsetsymbol{equal}{*}

\begin{icmlauthorlist}
\icmlauthor{Brian Irwin}{UBC}
\icmlauthor{Eldad Haber}{UBC}
\icmlauthor{Raviv Gal}{IBM}
\icmlauthor{Avi Ziv}{IBM}
\end{icmlauthorlist}

\icmlaffiliation{UBC}{Department of Earth, Ocean, and Atmospheric Sciences, University of British Columbia, Vancouver, 
Canada}
\icmlaffiliation{IBM}{Hybrid Cloud Quality Technologies, IBM Research, Haifa, Israel}

\icmlcorrespondingauthor{Brian Irwin}{birwin@eoas.ubc.ca}

\icmlkeywords{Neural Network Algorithms, Applications, Optimization}

\vskip 0.3in
]



\printAffiliationsAndNotice{}  

\begin{abstract}
In this paper, we illustrate a novel method for solving optimization problems when derivatives are not explicitly available. We show that combining implicit filtering (IF), an existing derivative free optimization (DFO) method, with a deep neural network global approximator leads to an accelerated DFO method. Derivative free optimization problems occur in a wide variety of applications, including simulation based optimization and the optimization of stochastic processes, and naturally arise when the objective function can be viewed as a black box, such as a computer simulation. We highlight the practical value of our method, which we call deep neural network accelerated implicit filtering (DNNAIF), by demonstrating its ability to help solve the coverage directed generation (CDG) problem. Solving the CDG problem is a key part of the design and verification process for new electronic circuits, including the chips that power modern servers and smartphones.
\end{abstract}

\section{Introduction}
\label{sec:intro}
Derivative free optimization (DFO), also known as zeroth order (ZO) optimization, is an important subfield within numerical optimization with many practical applications that range from fluid dynamics and finance to hardware verification. For examples, see \cite{DFONoisyFunctionsQuasiNewton}, \cite{DFOAlgReviewSoftwareComparison}, \cite{kelley3}, \cite{Gal_2020_HowToCatchALion} and the recent review in \cite{Larson2019DerivativefreeOM}. In such problems, the function that is optimized is often viewed as a black box, whose inner workings are unknown. Given an input, the black box produces an output, typically with some added noise. Furthermore, the nature of the noise may vary. For example, the noise can be systematic and deterministic, in which case the noise may be as simple as adding a high frequency periodic function to the output, resulting in samples from the black box producing a rough output surface. Alternatively, the noise may be stochastic, in which case the value of the noise changes each time we sample from the black box, even if the inputs are unchanged. As we do not know the inner workings of the black box, it becomes impossible to analytically compute derivatives. In addition, if the output of the black box is stochastic, it may be difficult, if not impossible, to use derivative based methods due to the randomness. 

In this paper, we consider such applications where derivative based optimization methods are not feasible. In particular, we focus on problems where output is corrupted by stochastic noise. While there are a variety of techniques for dealing with noisy optimization problems where derivatives are not available, we explore the use and acceleration of methods that are based on implicit filtering \cite{kelley3}. These methods are designed to work with noisy functions, and have been proven to give useful results in fields that range from control problems in ground water optimization \cite{GW2001} to hardware verification \cite{Gal_2020_HowToCatchALion}. 

Before briefly reviewing implicit filtering (IF), we identify key characteristics of the problems we aim to solve. First, we assume that one can sample the function, typically in parallel. Second, we assume that evaluating the function is very expensive. By expensive, we mean evaluating the function is resource intensive, which may mean it takes a very long time to evaluate the function, or may require a large amount of computational resources to compute the function, or both. Problems where function evaluation is expensive typically involve complicated simulations or large datasets, and often contain stochastic elements. 

We now briefly review the fundamentals of implicit filtering. Implicit filtering can be viewed as a trust region based direct search method. At each iteration of the algorithm, one searches for a better point by sampling the function along the boundary of a trust region. To illustrate, searching for a better point may consist of sampling the function along the vertices of a grid shaped trust region. When a better point is found, the algorithm takes that point, and a new trust region centred at the new point is used for future searches until an even better point is found. On the other hand, if the search along the boundary of the trust region fails to find a better point, the algorithm reduces the trust region size and tries again to look for a better point along the boundary of the smaller trust region. 

As implicit filtering only involves sampling the function directly without computing the gradient, it is attractive for problems where derivatives are unavailable. One obvious improvement that is thoroughly discussed in \cite{kelley3} involves using the sampled points at each iteration to build a local quadratic approximation of the function, and then attempting to perform a Newton-like update using the local quadratic approximation in order to accelerate convergence. This improvement can work well if the local quadratic approximation is sufficiently accurate. However, the local quadratic approximation is often not sufficiently accurate until relatively late iterations of implicit filtering, at which point the trust region is small and a quadratic approximation to a complex surface is sufficient.

A major disadvantage of implicit filtering is that the algorithm is memoryless. Indeed, each iteration of the implicit filtering algorithm ignores function evaluations from previous iterations. Assuming storage is relatively cheap compared to evaluating the function, not using the (typically many) function evaluations computed during the course of implicit filtering seems wasteful. Previously computed values of the function, even in places far from an optimum, typically contain useful information about the structure of the function. When used properly, such information can be harnessed to decrease the amount of computation required when optimizing the function.

In this work, we present an approach for leveraging previously computed function values to accelerate the progress of implicit filtering type methods. In particular, we use a deep neural network (DNN) to build an approximate model of the function using all the computed values of the function so far. The DNN approximation is then used to suggest new sampling points distinct from the implicit filtering trust region search, which can be much better than the points sampled by implicit filtering alone. We show that such an approach can at best dramatically improve over implicit filtering, while at worst having almost the same speed as standard implicit filtering.  

Building a deep neural network approximation to the function may be an expensive task on its own, depending on the amount of data required for the DNN approximation to become sufficiently accurate, and the amount of computations required when optimizing the DNN parameters so the DNN output fits the data. However, for the problems we focus on in this paper, fitting and evaluating the DNN approximation is cheap compared to evaluating the true function. In this case, using a DNN to build an approximation that can be used as a proxy or surrogate in place of the true function can ultimately reduce the number of expensive true function evaluations, and thus translate to a large computational gain overall. 

In addition, one may argue that it is possible to not use implicit filtering at all, and instead simply sample the function, approximate it using a DNN surrogate model, and minimize the DNN surrogate model. However, this approach is not competitive neither in theory nor in practice. As the DNN may converge to the sampled function values slowly, the overall convergence of this approach may be slow. Furthermore, to the authors' knowledge, there is no proof of convergence for such an approach. On the other hand, implicit filtering has a strong theoretical background, complete with established convergence properties (see \cite{kelley3}). Therefore, improving implicit filtering by including search points that are chosen by a DNN surrogate allows one to enjoy both worlds, taking points that are suggested by the DNN (exploiting the surrogate) whenever they provide sufficient progress, and relying on the points obtained by implicit filtering (exploring the space) when the DNN performs poorly. As evaluating the function is expensive, we emphasize that for the problems this paper focuses on, compared to other areas of machine learning, the number of sampled points is relatively small, and thus it is unreasonable to assume the DNN surrogate will yield a very good approximation to the function everywhere. 

The rest of this paper is organized as follows. First, we review implicit filtering and the deep neural networks used in this work. Next, we discuss how to use deep neural networks to accelerate the implicit filtering algorithm and discuss implementation issues and computational cost. After these discussions, we present a convergence result, and then perform numerical experiments. The first numerical experiment uses a small model problem that allows us to gain some fundamental understanding, and the second numerical experiment uses a realistic problem, the coverage directed generation (CDG) problem taken from the field of hardware verification. Finally, we summarize and conclude the paper.

\section{Implicit Filtering and Deep Neural Networks}
\label{sec:IF-and-DNNs}
In this section, we review the background for both implicit filtering and deep neural networks, and show how to combine the two approaches in order to obtain a robust optimization technique that leverages the advantages of both approaches.
Our goal is to solve the unconstrained  optimization problem
\begin{eqnarray}
\label{min}
\min_{\bfx} \big \{ f(\bfx) \big \} 
\end{eqnarray}
where $\bfx \in \mathbb{R}^n$ is a vector and $f(\cdot)$ is a function that typically contains noise and whose derivatives are infeasible to obtain. Bold uppercase letters (e.g. \bfI) represent matrices and bold lowercase letters (e.g. \bfx) represent vectors.

\subsection{Background - Implicit Filtering}
\label{sec:IF-background}
As previously mentioned, one attractive way to solve noisy optimization problems where derivative information is lacking or unavailable, such as in (\ref{min}), is implicit filtering. Below, we present a very simplified version of implicit filtering. For a thorough discussion, see \cite{kelley3} where the algorithm is introduced, applied to several problems, and its convergence properties are reviewed. Given a current point $\bfx_k$, a simple way to find a possible better point is to sample the function in a trust region around the current point. To this end, we evaluate the function at points $\bfx_k+h\bfw_j, j=1,\ldots, n_s$, where $\bfw_j$ are, for example, points on a grid or points drawn randomly from some predefined distribution, such as the Rademacher distribution or uniformly on the surface of a hypersphere. The parameter $h$ represents the trust region (or stencil) size and defines the resolution of the search, while $n_s$ is the number of search directions (i.e. new points to sample at). A search is considered successful if for one of the points, $j^*$,  we have that $f(\bfx_k+h\bfw_{j^*}) < f(\bfx_k)$. In this case, we set $\bfx_{k+1} = \bfx_k+h\bfw_{j^*}$ and continue. If none of the points yield a reduction in the function value, a stencil failure has occurred, and we set $h \leftarrow \tau_{tr} h$, where $\tau_{tr} \in (0,1)$, and repeat the process. The algorithm is summarized in Algorithm \ref{algPP}.
\begin{algorithm}[h]
\caption{ Implicit Filtering   \label{algPP}}
\begin{algorithmic}[1]
    \STATE Set $h=h_0$ and choose an initial point $\bfx_0$ and compute $f_0 = f(\bfx_0)$
    \WHILE{$h>h_{\rm min}$}
        \STATE Choose $n_s$ directions and set \\ $\bfX_k = [\bfx_k+h\bfw_{1},\ldots, \bfx_k+h\bfw_{n_s}]$
        \STATE Compute the values of the function \\ $\bff_k = [f(\bfx_k+h\bfw_{1}),\ldots,f( \bfx_k+h\bfw_{n_s})]^{\top}$
        \STATE Set $f_{\rm try} = \min \big \{ \bff_k \big \}$ and $j^*$ the index of $f_{\rm try}$
         \IF{$f_{\rm try} < f_k $}
         \STATE Set $\bfx_{k+1} = \bfx_k+h\bfw_{j^*}$
        \ELSE  
        \STATE Shrink $h$ upon stencil failure
        \STATE $h \leftarrow \tau_{tr} h$
        \ENDIF
    \ENDWHILE
\end{algorithmic}
\end{algorithm}

Note that at every implicit filtering iteration, we sample the function at points $\bfX_k = [\bfx_k+h\bfw_{1},\ldots, \bfx_k+h\bfw_{n_s}]$, and we know the values of the function at $\bff_k = [f(\bfx_k+h\bfw_{1}),\ldots,f( \bfx_k+h\bfw_{n_s})]^{\top}$. These values can be used to obtain a local approximation to the function $f$ (for example, by building a linear or quadratic approximation as discussed in \cite{kelley3}).  However, in the standard implicit filtering algorithm, previous points $[\bfX_1, \ldots, \bfX_{k-1}]$ and their function values $[\bff_1,\ldots,\bff_{k-1}]$ are not used in future calculations. If storage is available, such misuse is wasteful, as these values may hold valuable information about the function that can lead to more robust convergence. We now discuss how to leverage these values using a neural network.

\subsection{Background - Deep Neural Networks} 
\label{sec:DNN-background}
Our goal is to build a relatively inexpensive surrogate model of the function $f(\bfx)$ such that it can be probed to obtain an approximate minimizer of $f(\bfx)$. Neural networks, and in particular deep neural networks, are a class of function approximators that use decomposition to approximate a function. In this work, we chiefly use residual networks \cite{he2016identity,LiEtAl2017}, which have proven to be easy to train for difficult and highly nonlinear tasks. Our networks have the general structure
\begin{eqnarray}
\label{network}
\bfy_0 &=& \bfx \\
\nonumber
\bfy_1 &=& \bfG_0 \sigma_0 (\bfK_0 \bfy_0 + \bfb_0) \\
\nonumber
\bfy_{j+1} &=& \bfW_j\bfy_j + \bfG_j \sigma_{j} (\bfK_j \bfy_j + \bfb_j), \text{  } j=1,\ldots, N-2 \\
\nonumber
\bfy_N &=& \bfK_{N-1} \bfy_{N-1} + \bfb_{N-1}
\end{eqnarray}
Here, we choose $\bfW_j=\bfI$, and the dimensions of $\bfy_j$ and $\bfy_{j+1}$ are always consistent. The $N$ sets of parameters $\bftheta_j = \{\bfG_j, \bfK_j, \bfb_j \}, j=0,\ldots,N-2$ and $\bftheta_{N-1} = \{\bfK_{N-1}, \bfb_{N-1} \}$ are trainable parameters. For the networks employed in this work, the $\bfK_j$ are dense matrices. However, for other applications, one may consider sparse matrices, such as convolutions or other special structures. In our implementations, we chose either $\bfG_j = \alpha \bfI$ or $\bfG_j = -\alpha \bfK_j^{\top}$ with $\alpha > 0$. We chiefly set the activation functions $\sigma_j(\cdot)$ as the $\rm relu$ function. This type of network is stable and robust as shown in \cite{HaberRuthotto2017a}, especially when the intrinsic dimension of the output is smaller than the input. For brevity, we concatenate all the trainable network parameters into a single notation $\bfTheta = \{\bftheta_0,\ldots,\bftheta_{N-1}\}$. 

By choosing $\bfy_0 = \bfx$, we can forward propagate to obtain a vector $\bfy_N(\bfx,\bfTheta) \in \mathbb{R}^m$, where $m$ is the final dimension of the output. Using the vector $\bfy_N(\bfx,\bfTheta)$, we now want to obtain a surrogate $\widehat f(\bfx,\bfTheta)$ to the original function $f(\bfx)$. 
There are many different approaches to obtain such an approximation. For example, it is possible to choose a quadratic model,
$\widehat f(\bfx,\bfTheta) = \hf \| \bfy_N(\bfx,\bfTheta) \|_2^2$, 
which is bounded from below. One can also simply set $m = 1$ and choose $f(\bfx,\bfTheta) = \bfy_N(\bfx,\bfTheta)$. Choices of $f(\bfx,\bfTheta)$ that incorporate specific information about the true function $f(\bfx)$ can also be used. Once chosen, the surrogate $\widehat f(\bfx,\bfTheta)$ can then be used to probe the original function, and can be fit using a loss function $\Psi(\bfx,\bfTheta)$, such as 
\begin{eqnarray} \label{eq:squared-residual-loss}
\Psi(\bfx,\bfTheta) = \frac 12\, \bfr(\bfx,\bfTheta)^{\top} \bfr(\bfx,\bfTheta)
\end{eqnarray}
where $\bfr(\bfx,\bfTheta)$ is some residual. We can now use points $\bfX$ and measured function values $\bff$ in order to train the network and assess its parameters $\bfTheta$ by minimizing $\Psi(\bfx,\bfTheta)$. In this work, we use a mean squared error in order to estimate the network parameters. In~(\ref{eq:squared-residual-loss}), this corresponds to choosing $\bfr(\bfx,\bfTheta) = \bff(\bfx) - \widehat \bff(\bfx,\bfTheta)$. Once the network is trained, we have an approximation to the true function $f(\bfx)$ that can be used for the purpose of accelerating implicit filtering.

\section{Accelerating Implicit Filtering}
\label{sec:accelerating-IF}
Let us now combine implicit filtering with the DNN surrogate model. We start by using the originally obtained points $\bfX_0$ and $\bff_0$ to train a deep neural network. At every iteration of the implicit filtering algorithm, we obtain $n_s$ sampling points $\bfX_k$ and function values $\bff_k$. These points can now be added to the set of points previously calculated, and used to retrain the network by minimizing a loss function, giving an approximately optimal parameter vector $\bfTheta_k$. After training, we hold an approximate surrogate $\widehat f(\bfx,\bfTheta_k)$. This surrogate can be used in order to obtain an approximate minimizer of the function $f(\bfx)$. To this end, we use a gradient descent type technique to propose a new point. At each iteration $k$, we use a fixed number $s$ of iterations of the form
\begin{eqnarray}
\label{descent}
\bfx^{0}_k &=& \bfx_k \\
\nonumber
\bfx^{\ell+1}_k &=& \bfx_k^{\ell} - \mu_k\grad_{\bfx} \widehat f(\bfx_k^{\ell},\bfTheta_k)
\end{eqnarray}
with $\ell = 0, \dots, s-1$. The gradient of the surrogate with respect to the function argument $\bfx$ is computed using automatic differentiation (AD) \cite{nw}. The learning rate $\mu_k$ is computed using an Armijo line search \cite{nw}. During the iteration, we track the size of the update $\delta \bfx = \bfx^{\ell+1}_k - \bfx_k$. The iteration terminates if $\|\delta \bfx\| > h$, where $\| \cdot \|$ is an appropriate norm, which guarantees that we stay within the trust region of the implicit filtering method. While this may seem unimportant in practice, it is necessary to formally prove the convergence of the method (see Section~\ref{sec:convergence}). The discussion above is summarized in Algorithm \ref{dec}.
\begin{algorithm}[h]
\caption{ Surrogate Function Gradient Descent  \label{dec}}
\begin{algorithmic}[1]
    \STATE Choose $s$ and $h$, set $\ell=0$ and $\bfx^{0}_k = \bfx_k$
    \WHILE{$\ell < s$ and $ \|\delta \bfx \| \leq h$}
        \STATE Using AD, compute  ${\bf g}_k^{\ell} = \grad_{\bfx} \widehat f(\bfx_k^{\ell},\bfTheta_k)$
        \STATE Set $ \bfx^{\ell+1}_k = \bfx_k^{\ell} - \mu_k {\bf g}_k^{\ell}$ 
         \STATE Compute $\delta \bfx =  \bfx^{\ell+1}_k - \bfx_k$
         \STATE $\ell \leftarrow \ell+1$
        \ENDWHILE
\end{algorithmic}
\end{algorithm} 

The goal of the iteration in (\ref{descent}) is not to obtain the exact minimum of the DNN surrogate model $\widehat f(\bfx,\bfTheta_k)$, as the surrogate model is only an approximation of the true function $f(\bfx)$. A reduction in the true function value is sufficient. Therefore, we typically use a small number of steps $s$ and make sure that the surrogate value decreases. We then set $\bfx_{\rm try} = \bfx^{s}_k$. This point clearly reduces the surrogate value, and thus we use it as a potential candidate point for improving upon the value at the current point. This use of the surrogate is obvious, and it can improve upon a simple random or grid search. 

However, there is another important use of the surrogate model that can potentially improve the overall performance of the algorithm. An implicit filtering iteration requires choosing new search points. These points are chosen without using any form of filter or quality control. Without a filter or some valid assumption, the quality of the new search points is basically random. However, given a surrogate model, we can use the surrogate model to propose new points that are more likely to be of good quality and decrease the true function. Rather than choosing a point $\bfx_k + h\bfw_j$ and computing its value $f(\bfx_k + h\bfw_j)$, we first compute its surrogate value $\widehat f(\bfx_k + h\bfw_j, \bfTheta_k)$ and use this value as part of a filter. If
\begin{eqnarray}
\widehat f(\bfx_k + h\bfw_j, \bfTheta_k) \le f(\bfx_k) ,
\end{eqnarray}
then we pass the point to the true function and test the value of $f(\bfx_k + h\bfw_j)$. This is summarized in Algorithm~\ref{alg2}.
\begin{algorithm}[h]
\caption{ Surrogate Filtered Sampling  \label{alg2}}
\begin{algorithmic}[1]
    \STATE Choose the number of points needed $n_f$ and number of search directions $n_s$, set $\bfX_k = []$, $n=0$, and $j=0$
    \WHILE{$n<n_f$ and $j < {\rm n_s}$}
        \STATE Choose a direction $\bfw_j$
        \STATE Compute the value of the surrogate  $\widehat f(\bfx_k + h\bfw_j, \bfTheta_k)$
         \IF{$\widehat f(\bfx_k + h\bfw_j, \bfTheta_k) \le f(\bfx_k) $}
         \STATE Accept the point \\ $\bfX_{k}[:,n] = \bfx_k+h\bfw_j$
         \STATE $n \leftarrow n+1$
        \ENDIF
        \STATE $j \leftarrow j+1$
       \ENDWHILE
\end{algorithmic}
\end{algorithm}

Using the surrogate model to filter points can be very efficient when the network provides a reasonably good approximation to the true function. However, this may not always be the case, especially during the initial steps of the implicit filtering algorithm when the number of points used for training the surrogate DNN is small. Therefore, we propose to divide the number of sampling points into two groups. The first group, which we call the exploration group, chooses points on a grid or at random. The second group, which we call the exploitation group, uses the surrogate to choose sampling points. We have found experimentally that starting with a large number of exploration points, and then reducing them to about $20\%$ of the total points, yields good results. More work is needed in order to determine the best ratios of exploration to exploitation points as the algorithm progresses.

Finally, we combine the surrogate sampling and optimization together into an algorithmic framework we name Deep Neural Network Accelerated Implicit Filtering (DNNAIF). 
\begin{algorithm}[h]
\caption{ DNN Accelerated Implicit Filtering   \label{DNNAIF}}
\begin{algorithmic}[1]
    \STATE Set $h=h_0$ and choose an initial point $\bfx_0$ 
    \STATE Compute $f_0 = f(\bfx_0)$
    \STATE Train a network to obtain $\widehat f(\bfx, \bfTheta_0)$
    \WHILE{$h>h_{\rm min}$}
    	\STATE Use Algorithm \ref{dec} and propose a try point $\bfx_{\rm try}^0$
	\STATE Compute $f_{\rm try} = f( \bfx_{\rm try}^0)$
	\STATE Set $\bfX_k = \bfx_{\rm try}^0, \bff_k = f_{\rm try}$
	\IF{$f_{\rm try} < f_k$}
	\STATE $\bfx_{k+1} = \bfx_{\rm try}$
	\ELSE
        \STATE Choose $n_e$ directions and set \\ $\bfX_k^e = [\bfx_k+h\bfw_{1},\ldots, \bfx_k+h\bfw_{n_e}]$
        \STATE Choose $n_f$ points $\bfX_k^f$ using Algorithm~\ref{alg2}
        \STATE Set $\bfX_k = [\bfx_{\rm try}^0, \bfX_k^e, \bfX_k^f]$ 
        \STATE Compute the values of the function $\bff_k = f(\bfX_k)$
        \STATE Set $f_{\rm try} = \min \big \{ \bff_k \big \}$ and $j^*$ the index of $f_{\rm try}$
         \IF{$f_{\rm try} < f_k $}
         \STATE Set $\bfx_{k+1} = \bfX_k[:,j^*]$
        \ELSE  
        \STATE Shrink $h$ upon stencil failure
        \STATE $h \leftarrow \tau_{tr} h$
        \ENDIF
        \ENDIF
        \STATE Use previous data $\bfX_j, \text{  } j = 0,\dots,k-1$ and/or current data $\bfX_k$ to retrain the network and obtain $\widehat f(\bfx, \bfTheta_k)$
    \ENDWHILE
\end{algorithmic}
\end{algorithm}
The basic idea behind Algorithm~\ref{DNNAIF} is to use the surrogate function obtained by training a neural network as much as possible. However, since the amount of data that is needed for the network to give sufficiently accurate results may be rather large, we rely on the implicit filtering mechanism to compensate whenever the deep neural network approximation fails. We have found in most of our experiments that when the DNN surrogate is trained on a sufficient number of points, it significantly aids convergence initially, and then towards the end of the optimization process does not provide a significant advantage. In cases where a very accurate minima is desired, one should be encouraged to combine this approach with a local model, such as a local quadratic approximation to the function, as is often done with implicit filtering. Nonetheless, for noisy problems, where the accuracy needed may not be very high, early termination at a relatively large $h$ may be sufficient. In these cases, the algorithm can be very efficient and take many steps that are based on the surrogate. This can present significant overall savings compared to approaches that are based on the surrogate model or implicit filtering alone. 

Finally, an important point to consider is the overhead time that is added to the implicit filtering algorithm by using the surrogate model. This added time is due to training the network and approximately minimizing it. While this may not be a trivial amount of time, for the problems focused on in this paper where function evaluations are expensive and can require significant computational resources, this time is relatively insignificant, and well spent in the sense that it can significantly reduce the overall total computational time by avoiding true function evaluations. Furthermore, the cost of retraining the network can be reduced by using techniques such as warm starting and incremental training. 

\section{Convergence of DNNAIF}
\label{sec:convergence}
In this section, we outline how Algorithm~\ref{DNNAIF} inherits the asymptotic convergence properties of implicit filtering. For in-depth treatments of the convergence properties of implicit filtering, we refer the reader to \cite{kelley3}, \cite{Bortz1998}, and \cite{doi:10.1137/1.9781611970920}. We assume the objective function $f(\bfx)$ is bounded below and can be decomposed as 
\begin{equation} \label{eq:objective-decomposition}
f(\bfx) = f_s(\bfx) + \phi(\bfx)
\end{equation}
where $f_s(\bfx)$ is a smooth function and $\phi$ is noise. Following Chapter 5 of \cite{kelley3}, the key observation is that stencil failure (i.e. none of the function values on the stencil improve the base point) with a positive spanning set of directions is sufficient to conclude that $\norm{\nabla f_s}_2 = O(h)$. By positive spanning set (see  \cite{IntroductionToDFO}) of directions $\bfW = [ \bfw_1, \dots, \bfw_J ]$, we mean that any vector $\bfx \in \mathbb{R}^n$ can be written as 
\begin{equation}
\bfx = \sum_{j=1}^{J} a_j \bfw_j = \bfW \bfa 
\end{equation}
for some coefficient vector $\bfa \in \mathbb{R}^J$ where $a_j \geq 0, \forall j \in \{1,\dots,J\}$. Note that $\bfa$ is not necessarily unique. Denote the stencil $S$ defined by the set of directions $\bfW$ as 
\begin{equation} \label{eq:stencil-def}
S(\bfx,h,\bfW) = \{ \bfz \text{  } | \text{  } \bfz = \bfx + h \bfw_j, 1 \leq j \leq J \}
\end{equation}
and define the local norm of the noise as 
\begin{equation} \label{eq:noise-local-norm-def}
\norm{\phi}_{S(\bfx,h,\bfW)} = \max_{\bfz \in \{\bfx\} \cup S(\bfx,h,\bfW)} \big \{ \lvert \phi(\bfz) \rvert \big \}
\end{equation}
and the condition number of the positive spanning set $\bfW$ as
\begin{equation}
\kappa(\bfW) = \sqrt{n} \min \big \{ \norm{\bfA}_{\infty} \big \}
\end{equation}
where $\norm{\bfA}_{\infty} = \sup_{\bfx \neq 0} \big \{ \frac{\norm{\bfA \bfx}_{\infty}}{\norm{\bfx}_{\infty}} \big \}$ and the $2n \times J$ matrix $\bfA$ is constrained such that its entries are nonnegative and $\bfA \bfW^T = [\bfU, -\bfU]^T$ for some orthogonal matrix $\bfU$. With these definitions in hand, we now state a key theorem. The proof is given on pages 75 and 76 of \cite{kelley3}.

\begin{kelley-theorem*} \label{thm:kelley-theorem-5.7}
Let $f$ satisfy (\ref{eq:objective-decomposition}). Let $\nabla f_s$ be Lipschitz continuous with Lipschitz constant $L$. Let $\bfW$ be a positive spanning set. Then stencil failure implies that 
\begin{equation} \label{eq:stencil-failure-grad-bound}
\norm{\nabla f_s (\bfx)}_2 \leq \kappa (\bfW) \bigg ( \frac{L h}{2} + \frac{\norm{\phi}_{S(\bfx,h,\bfW)}}{h} \bigg )
\end{equation}
\end{kelley-theorem*}

As a result of (\ref{eq:stencil-failure-grad-bound}), Algorithm~\ref{DNNAIF} (i.e. DNNAIF) is guaranteed to converge to a critical point of $f_s(\bfx)$ when the conditions of Theorem 5.7 are satisfied, $h_{\rm min} = 0$, and 
\begin{equation} \label{eq:noise-decay}
\lim_{k \rightarrow \infty} \frac{\norm{\phi}_{S(\bfx_k,h_k,\bfW)}}{h_k} = 0  .
\end{equation}
Mathematically, 
\begin{equation} \label{eq:asymptotic-critical-point}
\liminf _{k \rightarrow \infty} \norm{\nabla f_s(\bfx_k)}_2 = 0  .
\end{equation}
To see this, observe that an iteration of Algorithm~\ref{DNNAIF} ends with either a decrease in the objective $f$ or stencil failure, and let $\{ \bfx_{k_i} \}$ be an infinite subsequence of $\{ \bfx_{k} \}$ for which the iteration terminates with stencil failure. As $h_{\rm min} = 0$, $h_{k_i} \rightarrow 0$ (i.e. $h$ eventually tends to zero), which combines with (\ref{eq:stencil-failure-grad-bound}) and (\ref{eq:noise-decay}) to give (\ref{eq:asymptotic-critical-point}).

\section{Numerical Experiments}
\label{sec:numerical-experiments}
In this section, we experiment with two different problems. The first problem is a simple model problem that we can easily explore and visualize. The second problem is a realistic problem that is taken from the field of hardware verification \cite{Gal_2020_HowToCatchALion}, and motivated us to develop the techniques discussed in this paper.

\subsection{The Noisy Rosenbrock Banana Function}
\label{sec:noisy-banana}
\begin{figure}
  \centering
  \begin{tabular}{cc}
  \includegraphics[width=3.75cm]{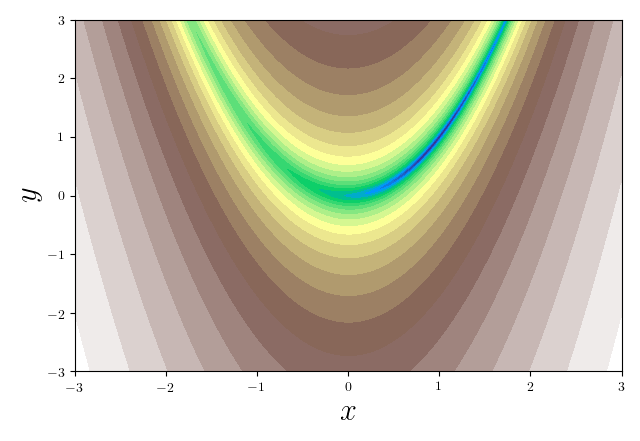} &
  \includegraphics[width=3.75cm]{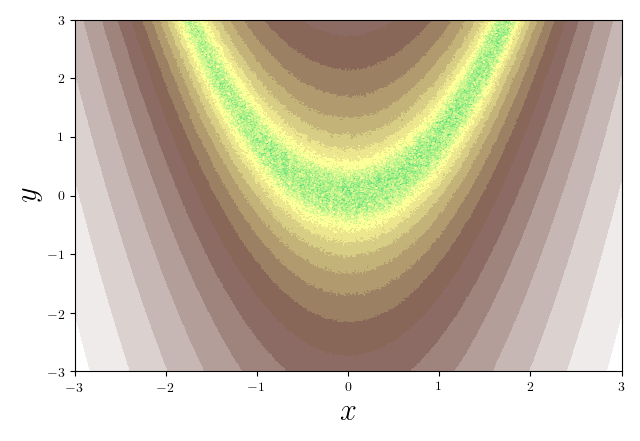} \\
  \end{tabular} 
  \caption{The Rosenbrock banana function (left) and a version with standard normal random noise added (right) used for numerical experiments. Observe how the noise obscures the location of the minimum in the green valley. Same colormap as Figures~\ref{fig:DNN_surrogate_init} and~\ref{fig:DNN_surrogate_fin_DNNAIF}. }
  \label{fig:banana}
\end{figure}

Minimizing the Rosenbrock banana function \cite{10.1093/comjnl/3.3.175} that maps the vector $\bfx = [x,y]^{\top}$ to $\mathbb{R}$ can be written in the form
\begin{equation}
\label{banana}
\min_{\bfx} \bigg \{ \bfr(\bfx)^{\top} \bfr(\bfx) \bigg \}
\end{equation}
where
\begin{equation}
\bfr(\bfx) = \begin{pmatrix} a - x\\ \sqrt{b} \big ( y - x^2 \big )  \end{pmatrix} .
\end{equation}
This problem can be thought of as a nonlinear data fitting problem. A noisy version of the problem can be written where we add independent and identically distributed (i.i.d.) Gaussian stochastic noise. The advantage of working with this problem is that we can visualize the behaviour of the algorithms and approximations. An image of the function and a version with standard normal random noise added is shown in Figure~\ref{fig:banana}.

\begin{figure}[h]
  \centering
  \includegraphics[width=9cm]{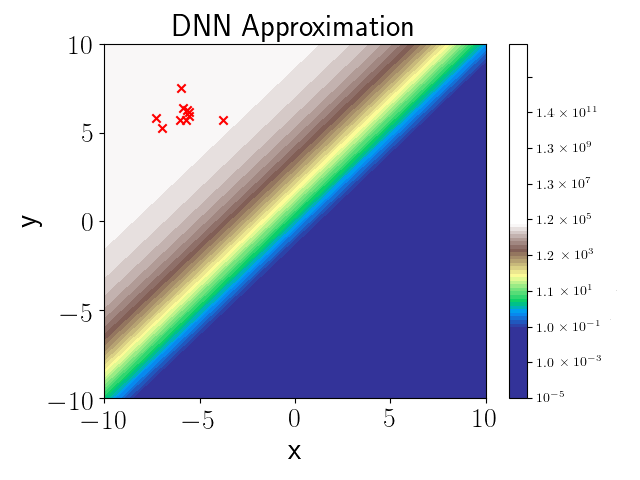} 
  \caption{Visualizing initial DNN approximation fit to randomly sampled initial data (red x marks). The fit roughly resembles a linear model, which can provide a reasonable descent direction. }
  \label{fig:DNN_surrogate_init}
\end{figure}

\begin{figure}[h]
  \centering
  \includegraphics[width=7cm]{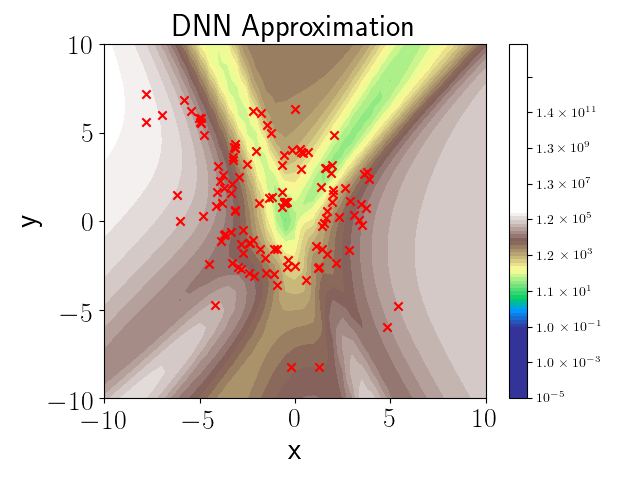} 
  \caption{Visualizing the final DNN approximation after 10 iterations of DNNAIF. The red x marks show data sampled during iterations 1 to 10. Observe how the final DNN approximation captures the general structure of the banana function. }
  \label{fig:DNN_surrogate_fin_DNNAIF}
\end{figure}

Before comparing the aggregate performance of IF, optimizing a DNN surrogate, and the combined approach of DNNAIF, we investigate and visualize the behaviour of the DNN surrogate models. In all cases, the DNN surrogate consists of a residual network with $30$ layers, and the network is trained using mini-batch gradient descent for $1000$ iterations. Figures~\ref{fig:DNN_surrogate_init} and~\ref{fig:DNN_surrogate_fin_DNNAIF} visualize the progress of the DNN surrogate for DNNAIF. Starting from a few random points sampled near $(-6,6)$, the initial DNN surrogate roughly resembles a linear model that does not capture the structure of the banana function well. However, this is not surprising given the small number of initial points. In contrast, as shown in Figure~\ref{fig:DNN_surrogate_fin_DNNAIF}, when combined with IF, the DNN is able to use the implicit filtering points as additional function exploration, which enables the DNN surrogate to more quickly become a useful model of the general structure of the function. 

\begin{figure}[h]
  \centering
  \includegraphics[width=7cm]{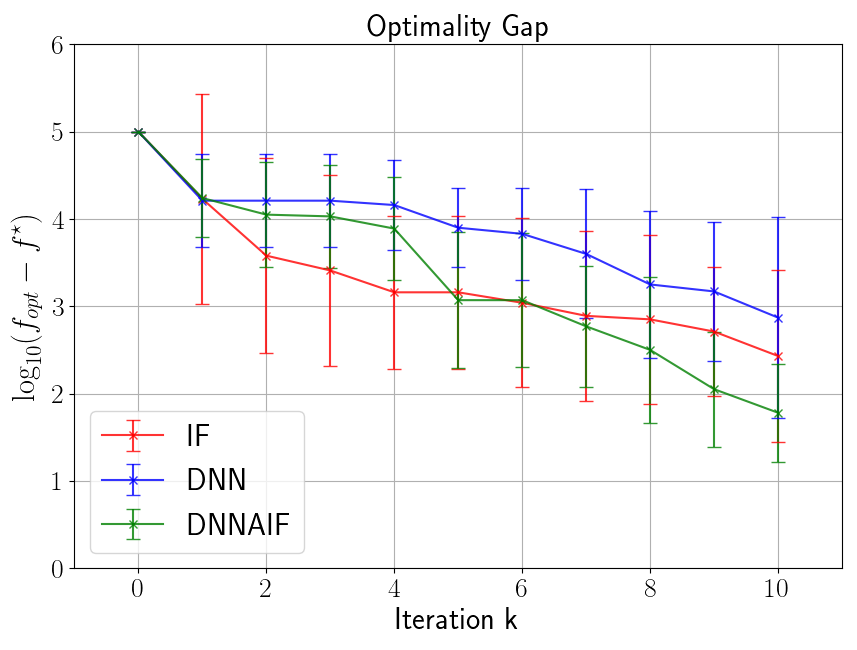} 
  \caption{Algorithm performance comparison. Starting from the point $(-6,6)$, each algorithm runs for 10 iterations, using 11 points at each iteration, and optimizes the noisy banana function. The metric plotted is the average optimality gap for the true (i.e. no noise) banana function in log space, with the error bars showing a single standard deviation. Statistics are taken over $10$ runs.  $f_{opt}$ is the true value of the best $f$ at each iteration, and $f^{\star} = 0$ is the minimum value of the non-noisy banana function. }
  \label{fig:opt_gap}
\end{figure}

Having visualized how the DNN surrogates can evolve and how the exploratory behaviour of IF can enhance the DNN approximation, we are now ready to compare the performance of the three approaches on the banana function. Figure~\ref{fig:opt_gap} uses the optimality gap on the true (i.e. non-noisy) banana function as a metric. For these experiments, IF and DNNAIF used a circular trust region with $h_0 = 30$ and $\tau_{tr} = 0.9$. No trust region was enforced for the DNN method. For DNNAIF, $10$ implicit filtering (exploration) points were used per iteration, with the $11^{th}$ point being a DNN surrogate suggested (exploitation) point. Here, we observe that both IF and DNNAIF outperform optimizing a DNN surrogate. We also observe that with this small number of sampling points, the benefits of DNNAIF do not become apparent until after iteration 6 or so. However, after this point, the surrogate model has a sufficient number of training points to start to capture the general structure well, and accelerates the algorithm so that the final average optimality gap for DNNAIF is smaller than the optimality gap for IF after $10$ iterations. Furthermore, the worst case performance of DNNAIF is superior to the worst case performance of IF at later iterations. 

When analyzing the results of this experiment, it is important to realize that IF can do very well due to the low dimensionality of the problem. As the dimension of the problem increases, IF requires more directions to have similar behaviour. As a result, the DNN acceleration becomes more crucial. This is explored in the next experiment.

\subsection{Application to Coverage Directed Generation}
\label{sec:CDG}
Coverage Directed Generation (CDG)~\cite{cdg, gilly_cdg} is a generic name used for a multitude of techniques and algorithms that are used in hardware verification to create tests for hitting so called coverage events \cite{cov_book}. While most approaches rely on random sampling, some AI algorithms and optimization techniques have been proposed for the solution of the CDG problem~\cite{cdg-bayesian-dac, 10.1145/3061639.3062324}. The problem with these algorithms is that they frequently require a prohibitively large number of simulations to explore the solution space. As simulators of today's advanced processors often take minutes to hours to perform a simulation, these approaches can simply take too long. Our attempt here is to make such algorithms more practical by reducing the number of simulations required, and thus the total time required.

\begin{figure}
  \centering
  \includegraphics[width=7cm]{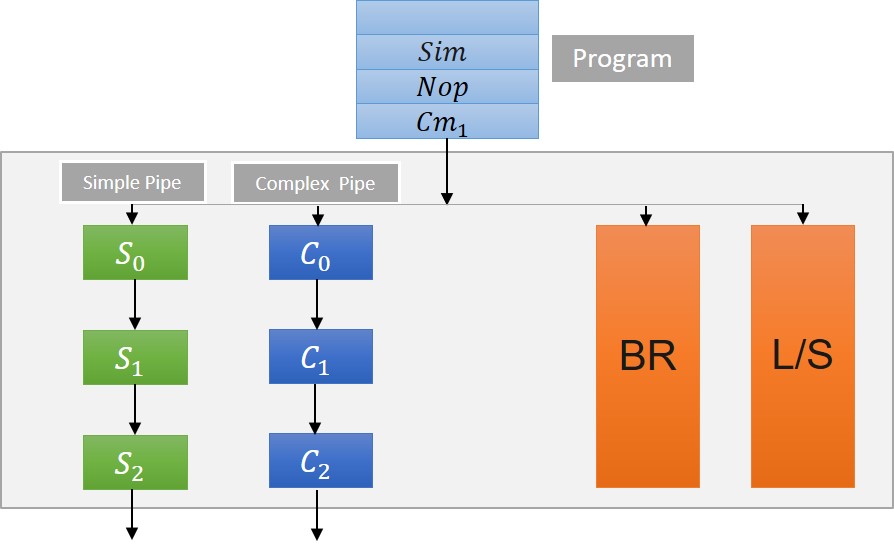}
  \caption{Schematic of the simulated NorthStar pipeline. There are two pipes of 3 stages: one simple pipe $S$ and one complex pipe $C$. In addition, L/S represents the processor's load store unit, and BR the branch prediction unit.}
  \label{fig:NS_HL}
\end{figure}

As an experimental environment, we employ an abstract high-level simulator of the two arithmetic pipes in the NorthStar \cite{Salvatore_4thgeneration} in-order processor and the dispatch unit, also used in \cite{cdg-bayesian-dac}. A sketch of the NorthStar pipeline is provided in Figure~\ref{fig:NS_HL}.

The goal of the verification process is to sufficiently sample every possible state or event of the system. Let $\bfx$ be the input to the NorthStar simulator. We follow \cite{Gal_2020_HowToCatchALion} and define the function $p_j(\bfx)$ as the probability of hitting the $j^{th}$-event. The goal is to maximize the probability by changing the inputs to the simulator, $\bfx$. We define the objective function to be maximized as  
\begin{equation}
f(\bfx) = \sum_{j} w(p_j(\bfx)) p_j(\bfx) 
\end{equation}
where $0\le w(t)$ is a monotonically decreasing weighting function. The goal of the weighting function $w(t)$ is to increase the reward for hitting unhit or hard-to-hit events and reduce the effect of easily hit events. For the CGD problem, the goal is to find simulation inputs that result in each event being hit with probability $p_j$ above a set threshold. Thus, the objective measures the cardinality of all events that are unhit.

Similar to the previous example, the objective function is noisy. The results of each simulation are random, as the input to the simulation defines a probability space rather than a fixed value that the system uses. To solve this problem, one does not solve a single optimization problem, but rather solves a sequence of optimization problems with the goal of covering all hard-to-hit events. We now explore the use of DNNAIF for the solution of the CDG problem with the NorthStar simulator, which has $23$ input parameters to optimize, and about $35$ events that are considered hard-to-hit. 

\begin{figure}[h]
  \centering
  \includegraphics[width=7cm]{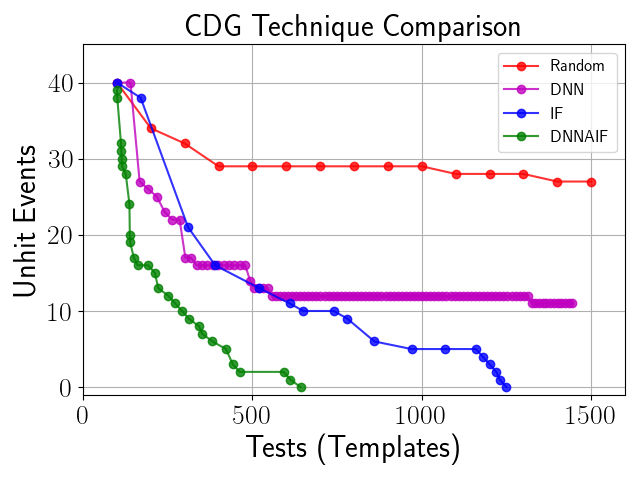}
  \caption{Comparison of four techniques for maximizing the events hit in the NorthStar simulator. The number of tests is equivalent to the number of true function evaluations, and the goal here is to hit all hard-to-hit events. $20$ sampling directions were used per iteration. For DNNAIF, $10$ directions were chosen at random and $10$ based on the DNN. Curves that hit zero faster are better. }
  \label{fig:NS_Res}
\end{figure}

Perhaps the most common approach in the CDG field is to use smart random sampling. For this particular problem, this often means sampling using a Dirichlet distribution. We compare this common approach of Dirichlet random sampling with the three optimization techniques from the banana function experiment, namely IF, optimizing a DNN surrogate, and the combined approach of DNNAIF. The results of these experiments are presented in Figure~\ref{fig:NS_Res}.
 
The results shown in Figure~\ref{fig:NS_Res} demonstrate that neither Dirichlet random sampling nor a DNN approximation to the function were the best way to solve the CDG problem, as many hard-to-hit events were never hit. IF and DNNAIF both hit all hard-to-hit events, but the number of true function evaluations for IF is rather large, and DNNAIF uses approximately half as many (about $600$ tests vs. $1200$ tests). In the world of hardware verification, where each true function evaluation may require massive computational resources, such savings in terms of true function evaluations is significant.

\section{Conclusions}
\label{sec:conclusion}
In this work, we have introduced a new derivative free optimization method, which we call DNNAIF, that combines a direct search method (implicit filtering) and an approximate surrogate function that is based on a deep neural network. DNNAIF can use the DNN surrogate to both find an acceptable new point at each iteration, and to aid the implicit filtering process by screening potential points before using them for true function evaluations. Both these uses help reduce the number of true function evaluations. Experiments with both the Rosenbrock function and NorthStar CDG problem establish that DNNAIF can provide superior performance in terms of true function evaluations compared to IF, as the authors intended when conceiving DNNAIF. 

Using DNNAIF requires choosing several parameters. In particular, the balance between exploration and exploitation can be problem dependent. For problems where the surrogate model trains quickly, one may use only a few iterations of implicit filtering with random sampling followed by swiftly using more points that are based on the surrogate model. The computational cost of DNNAIF at each iteration clearly depends on the choice of DNN architecture, the loss function, and the specific retraining procedure. This dependency comes from the need to update the DNN surrogate model at each iteration, unlike IF where there is no surrogate model to directly update. As a result, if one chooses poorly in these three areas, DNNAIF may not be competitive compared to other methods. Fortunately, the theoretical convergence guarantee for DNNAIF presented in Section~\ref{sec:convergence} does not depend on these three areas. Nonetheless, for many DFO problems of interest to the authors, such as the CDG problem and others in engineering design, even a small DNN may yield significant savings in terms of true function evaluations. In future work, we intend to address these issues by investigating the performance of DNNAIF using a wider variety of applications than can be considered in an introductory paper.  

\newpage

\bibliography{biblioFile}

\begin{thebibliography}{20}
\providecommand{\natexlab}[1]{#1}
\providecommand{\url}[1]{\texttt{#1}}
\expandafter\ifx\csname urlstyle\endcsname\relax
  \providecommand{\doi}[1]{doi: #1}\else
  \providecommand{\doi}{doi: \begingroup \urlstyle{rm}\Url}\fi

\bibitem[Battermann et~al.(2001)Battermann, Gablonsky, Patrick, Kelley,
  Kavanagh, Coffey, and Miller]{GW2001}
Battermann, A., Gablonsky, J., Patrick, A., Kelley, C., Kavanagh, K., Coffey,
  T., and Miller, C.
\newblock Solution of a groundwater control problem with implicit filtering.
\newblock \emph{Optimization and Engineering}, 01 2001.

\bibitem[Berahas et~al.(2019)Berahas, Byrd, and
  Nocedal]{DFONoisyFunctionsQuasiNewton}
Berahas, A.~S., Byrd, R.~H., and Nocedal, J.
\newblock Derivative-free optimization of noisy functions via quasi-newton
  methods.
\newblock \emph{SIAM Journal on Optimization}, 29:\penalty0 965--993, 2019.

\bibitem[Borkenhagen \& Storino(1999)Borkenhagen and
  Storino]{Salvatore_4thgeneration}
Borkenhagen, J. and Storino, S.
\newblock 4th generation 64-bit powerpc-compatible commercial processor design,
  1999.

\bibitem[Bortz \& Kelley(1998)Bortz and Kelley]{Bortz1998}
Bortz, D.~M. and Kelley, C.~T.
\newblock \emph{The Simplex Gradient and Noisy Optimization Problems}, pp.\
  77--90.
\newblock Birkh{\"a}user Boston, Boston, MA, 1998.
\newblock ISBN 978-1-4612-1780-0.
\newblock \doi{10.1007/978-1-4612-1780-0_5}.
\newblock URL \url{https://doi.org/10.1007/978-1-4612-1780-0_5}.

\bibitem[Conn et~al.(2009)Conn, Scheinberg, and Vicente]{IntroductionToDFO}
Conn, A., Scheinberg, K., and Vicente, L.
\newblock \emph{Introduction to Derivative-Free Optimization}.
\newblock SIAM, Philadelphia, 2009.

\bibitem[Fine \& Ziv(2003)Fine and Ziv]{cdg-bayesian-dac}
Fine, S. and Ziv, A.
\newblock Coverage directed test generation for functional verification using
  {B}ayesian networks.
\newblock In \emph{Proceedings of the 40th Design Automation Conference}, pp.\
  286--291, 2003.

\bibitem[Gal et~al.(2017)Gal, Kermany, Saleh, Ziv, Behm, and
  Hickerson]{10.1145/3061639.3062324}
Gal, R., Kermany, E., Saleh, B., Ziv, A., Behm, M., and Hickerson, B.
\newblock Template aware coverage: Taking coverage analysis to the next level.
\newblock In \emph{Proceedings of the 54th Annual Design Automation Conference
  2017}, DAC '17, New York, NY, USA, 2017. Association for Computing Machinery.
\newblock ISBN 9781450349277.
\newblock \doi{10.1145/3061639.3062324}.
\newblock URL \url{https://doi.org/10.1145/3061639.3062324}.

\bibitem[Gal et~al.(2020)Gal, Haber, Irwin, Saleh, and
  Ziv]{Gal_2020_HowToCatchALion}
Gal, R., Haber, E., Irwin, B., Saleh, B., and Ziv, A.
\newblock How to catch a lion in the desert: on the solution of the coverage
  directed generation ({CDG}) problem.
\newblock \emph{Optimization and Engineering}, May 2020.
\newblock \doi{10.1007/s11081-020-09507-w}.
\newblock URL \url{https://doi.org/10.1007%2Fs11081-020-09507-w}.

\bibitem[Haber \& Ruthotto(2017)Haber and Ruthotto]{HaberRuthotto2017a}
Haber, E. and Ruthotto, L.
\newblock Stable architectures for deep neural networks.
\newblock \emph{Inverse Problems}, 34\penalty0 (1), 2017.

\bibitem[He et~al.(2016)He, Zhang, Ren, and Sun]{he2016identity}
He, K., Zhang, X., Ren, S., and Sun, J.
\newblock Identity mappings in deep residual networks.
\newblock In \emph{European Conference on Computer Vision}, pp.\  630--645.
  Springer, 2016.

\bibitem[Kelley(2011)]{kelley3}
Kelley, C.
\newblock \emph{Implicit Filtering}.
\newblock SIAM, Philadelphia, 2011.

\bibitem[Kelley(1999)]{doi:10.1137/1.9781611970920}
Kelley, C.~T.
\newblock \emph{Iterative Methods for Optimization}.
\newblock Society for Industrial and Applied Mathematics, 1999.
\newblock \doi{10.1137/1.9781611970920}.
\newblock URL \url{https://epubs.siam.org/doi/abs/10.1137/1.9781611970920}.

\bibitem[Larson et~al.(2019)Larson, Menickelly, and
  Wild]{Larson2019DerivativefreeOM}
Larson, J., Menickelly, M., and Wild, S.~M.
\newblock Derivative-free optimization methods.
\newblock \emph{Acta Numer.}, 28:\penalty0 287--404, 2019.

\bibitem[Li et~al.(2018)Li, Xu, Taylor, Studer, and Goldstein]{LiEtAl2017}
Li, H., Xu, Z., Taylor, G., Studer, C., and Goldstein, T.
\newblock Visualizing the loss landscape of neural nets.
\newblock In \emph{Advances in Neural Information Processing Systems}, pp.\
  6389--6399, 2018.

\bibitem[Nativ et~al.(2001)Nativ, Mittermaier, Ur, and Ziv]{gilly_cdg}
Nativ, G., Mittermaier, S., Ur, S., and Ziv, A.
\newblock Cost evaluation of coverage directed test generation for the ibm
  mainframe.
\newblock In \emph{Proceedings of the 2001 International Test Conference}, pp.\
   793--802, October 2001.

\bibitem[Nocedal \& Wright(1999)Nocedal and Wright]{nw}
Nocedal, J. and Wright, S.
\newblock \emph{Numerical Optimization}.
\newblock Springer, New York, 1999.

\bibitem[Piziali(2004)]{cov_book}
Piziali, A.
\newblock \emph{Functional Verification Coverage Measurement and Analysis}.
\newblock Springer, 2004.

\bibitem[Rios \& Sahinidis(2013)Rios and
  Sahinidis]{DFOAlgReviewSoftwareComparison}
Rios, L.~M. and Sahinidis, N.~V.
\newblock Derivative-free optimization: a review of algorithms and comparison
  of software implementations.
\newblock \emph{Journal of Global Optimization}, 56:\penalty0 1247--1293, 2013.

\bibitem[Rosenbrock(1960)]{10.1093/comjnl/3.3.175}
Rosenbrock, H.~H.
\newblock {An Automatic Method for Finding the Greatest or Least Value of a
  Function}.
\newblock \emph{The Computer Journal}, 3\penalty0 (3):\penalty0 175--184, 01
  1960.
\newblock ISSN 0010-4620.
\newblock \doi{10.1093/comjnl/3.3.175}.
\newblock URL \url{https://doi.org/10.1093/comjnl/3.3.175}.

\bibitem[Ur \& Yadin(1999)Ur and Yadin]{cdg}
Ur, S. and Yadin, Y.
\newblock Micro-architecture coverage directed generation of test programs.
\newblock In \emph{Proceedings of the 36th Design Automation Conference}, pp.\
  175--180, June 1999.

\end{thebibliography}
\bibliographystyle{icml2021}

\end{document}